\documentclass[11pt, final]{article}
\usepackage{a4}
\usepackage{amsmath}%
\usepackage{amstext}%
\usepackage{amssymb}%
\usepackage{showkeys}%
\usepackage{epsfig}%
\usepackage{cite}
\newtheorem{theorem}{Theorem}

\newtheorem{lemma}[theorem]{Lemma}

\newtheorem{remark}{Remark}

\newcounter{unnumber}

\newcommand{\R}{\mathbb{R}}%
\newcommand{\N}{\mathbb{N}}%
\DeclareMathOperator*\inte{int}%
\DeclareMathOperator*\ri{ri}%
\DeclareMathOperator*\cl{cl}%
\DeclareMathOperator*\epi{epi}%
\DeclareMathOperator*\dom{dom}%

\DeclareMathOperator*\B{\overline{\R}}%
\DeclareMathOperator*\diam{diam}%

\title{Error bound results for convex inequality systems via conjugate duality}

 \author{Radu Ioan Bo\c{t} \thanks
 {Faculty of Mathematics, Chemnitz University of Technology,
D-09107 Chemnitz, Germany, e-mail:
 radu.bot@mathematik.tu-chemnitz.de. Research partially supported by DFG (German Research Foundation), project WA 922/1-3.} \and Ern\"{o} Robert Csetnek
 \thanks {Faculty of Mathematics, Chemnitz University of Technology,
D-09107 Chemnitz, Germany, e-mail:
 robert.csetnek@mathematik.tu-chemnitz.de}}

\begin{document}
\maketitle

\begin{center}\bf{Dedicated to Professor Marco A. L\' opez on the occasion
of his 60th birthday}\end{center}\vspace{0.5cm}

\noindent \textbf{Abstract.} The aim of this paper is to implement some new techniques, based on conjugate duality in convex optimization, for proving the existence of global error bounds for convex inequality systems. We deal first of all with systems described via one convex inequality and
extend the achieved results, by making use of a celebrated scalarization function, to convex inequality systems expressed by means of a general vector function. We also propose a second approach for guaranteeing the existence of global error bounds of the latter, which meanwhile sharpens the classical result of Robinson.\\



\section{Introduction and preliminaries}

Consider $(X,\|\cdot\|)$ a real normed space and $f:X\rightarrow \B=\R \cup \{\pm \infty\}$ a proper and convex function such that $S=\{x\in X:f(x)\leq 0\}$ is nonempty. We say that the \emph{global error bound} holds for the inequality \begin{equation}\label{ineq-sc}f(x)\leq 0, \ x\in X\end{equation} if there exists a constant $\alpha>0$ (depending on the initial data) such that \begin{equation}\label{err-b-sc}d(x,S)\leq \alpha f(x)_+ \ \forall x\in X\end{equation} where, for $\gamma\in\B$, $\gamma_+=\max\{\gamma,0\}$.

In this article we first implement some new techniques for proving the existence of a global error bound for \eqref{ineq-sc} under classical assumptions by making use of the convex conjugate duality.

Then we consider for a further real normed space $(Y,\|\cdot\|)$,  partially ordered by the nonempty convex closed cone $K \subseteq Y$, and $g :X \rightarrow Y^\bullet = Y \cup \{\infty_K\}$, a proper and $K$-convex function such that $Q=\{x \in X : g(x) \in -K\}$ is nonempty, the inequality system
\begin{equation}\label{ineq-vect}g(x)\leq_K 0, \ x\in X,\end{equation}
for which we say that the \textit{global error bound} holds if there exists a constant $\alpha>0$ (depending on the initial data) such that \begin{equation}\label{err-b-vect}d(x,Q)\leq \alpha d(g(x),-K) \ \forall x\in X.\end{equation}
The issue of the existence of error bounds for \eqref{ineq-vect}, a topic which has its roots in the paper of Hoffman \cite{hoffm} (see, also, \cite{lew-pang, klatte, pang, Zal-baikal}),  was investigated in the seminal paper of Robinson \cite{robin} and, by particularizing the statements given there, one can easily obtain corresponding assertions for the existence of global error bounds for \eqref{ineq-sc} (for $Y=\R$, $K=\R_+$ and $Y^\bullet = \R \cup \{+\infty\}$). In Section 3 we go into the opposite direction, namely, when having the global error bound results for \eqref{ineq-sc}, we show how one can derive corresponding statements for \eqref{ineq-vect} and make to this end use of an appropriate scalarization function. An alternative approach for guaranteeing the existence of a global error bound for \eqref{ineq-vect}, based on a result due to Simons \cite{simons}, is proposed in Section 4 and this provides a sharpening of the result of Robinson from \cite{robin}. We close the paper by deriving some conclusions and by proposing some topics for further research.

To make the paper self-consistent, we consider in the following some preliminary notions and results (see \cite{b-hab, bgw-book, EkTem, hiriart-lem1, Rock-carte, Zal-carte}). Having a real normed space  $(X,\|\cdot\|)$, we denote by $(X^*, \|\cdot\|_*)$ its topological dual space. By $\langle x^*, x\rangle=x^*(x)$ we denote the value of the continuous linear functional $x^*\in X^*$ at $x\in X$, while $B(0,1) = \{x \in X : \|x\| \leq 1\}$ and $B_*(0,1) = \{x^* \in X^* : \|x^*\|_* \leq 1\}$ are the \textit{closed unit balls} of $X$ and $X^*$, respectively. Given a subset $S$ of $X$, by $\inte S$, $\ri S$ and $\cl S$ we denote its \textit{interior}, \textit{relative interior} and \textit{closure}, respectively. The function $\delta_S : X \rightarrow \overline\R$, defined by $\delta_S(x) = 0$ for $x \in S$ and $\delta_S(x) = +\infty$, otherwise, is the \textit{indicator function} of $S$, while $\sigma_S : X^* \rightarrow  \overline\R$, defined by $\sigma_S(x^*) = \sup_{x \in S} \langle x^*,x \rangle$, is the \textit{support function} of $S$. Further, $d(\cdot, S) : X \rightarrow \overline{\R}$, $d(x,S) = \inf_{s \in S} \|x-s\|$, is the \textit{distance function} of the set $S$, and it is always Lipschitz continuous, being convex when $S$ is a convex set.

On $\overline{\R}$ we consider the following conventions: $(+\infty) - (+\infty) = +\infty$, $0(+\infty) = +\infty$ and $0(-\infty) = 0$. Having a function $f:X\rightarrow \overline\R$ we use the classical notations for its \textit{domain} $\dom f=\{x\in X:
f(x)<+\infty\}$, its \textit{epigraph} $\epi f=\{(x, r)\in X\times \R: f(x)\leq r\}$ and its \textit{lower level set at level $r \in \R$},  $L(f,r) = \{x \in X: f(x) \leq r\}$. The \textit{lower semicontinuous hull} of $f:X\rightarrow \overline\R$ is the function $\cl f:X\rightarrow \overline\R$ which has as epigraph $\cl(\epi f)$. We call $f$ \textit{proper} if $f(x)>-\infty$ for all $x\in X$ and $\dom f \neq\emptyset$. Further, by $f_+ : X\rightarrow \overline\R$ we denote the function defined by $f_+(x) =f(x)_+$ for $x \in X$. The \textit{conjugate function} of $f$ is $f^*:X^*\rightarrow \overline\R$ defined by $f^*(x^*)=\sup\{\langle x^*, x\rangle - f(x) : x\in X\}$, while the \textit{biconjugate function} of $f$ is $f^{**}:X^{**}\rightarrow \overline\R$ defined by $f^{**}(x^{**})=\sup\{\langle x^{**}, x^*\rangle - f^*(x^*) : x^*\in X^*\}$. For all $x \in X$ one has $f^{**}(x)=\sup\{\langle x^{*}, x \rangle - f^*(x^*) : x^*\in X^*\}$. Regarding a function and its conjugate we have the \textit{Young-Fenchel inequality} $f^*(x^*)+f(x)\geq \langle x^*, x\rangle$ for all $x\in X$ and $x^*\in X^*$.  When $S \subseteq X$, one obviously has $\delta_S^* = \sigma_S$.

When $f$ is a proper, convex and lower semicontinuous function, according to the \textit{Fenchel-Moreau Theorem}, one has $f(x) = f^{**}(x)$ for all $x \in X$. Given the proper functions $f,g :X\rightarrow \overline \R$, their \textit{infimal convolution} is the function $f \square g :X\rightarrow
\overline\R$, $(f \square g)(x)=\inf \{f(x-y) + g(y) : y \in X\}$. One has $(f \square g)^* = f^* + g^*$.

Having $(Y,\|\cdot\|)$ another real normed space, we call a set $K \subseteq Y$ \textit{cone} if for all $\lambda \geq 0$ and all $k \in K$ one has $\lambda k \in K$. For a given cone $K \subseteq Y$ we denote by $K^*=\{\lambda \in Y^*:\langle \lambda, k\rangle \geq 0\ \forall k\in K\}$ its \textit{dual cone}. A nonempty convex cone $K \subseteq Y$ induces on $Y$ a partial order ``$\leq_K$'', defined by $y \leq_K z \Leftrightarrow z-y \in K$ for $y,z \in Y$. To $Y$ we attach a greatest element with respect to ``$\leq_K$'', which does not belong to $Y$, denoted by $\infty_K$ and let $Y^\bullet:=Y\cup \{\infty_K\}$. Then for any $y\in Y^\bullet$ one has $y \leq_K \infty_K$ and we consider on $Y^\bullet$ the operations $y+\infty_K=\infty_K+y=\infty_K$ for all $y\in Y$ and
$t\cdot \infty_K=\infty_K$ for all $t\geq 0$. By convention, for every set $S \subseteq Y$, $d(\infty_K, S) := +\infty$ and  $\langle \lambda, \infty_K \rangle := +\infty$ for all $\lambda \in C^*$.

A function $h : Y^\bullet \rightarrow \overline{\R}$ is said to be \textit{$K$-increasing ($K$-decreasing)}, whenever for all $y,z \in Y$  with $y \leq_K z$ one has $h(y) \leq h(z)$ ($h(y) \geq h(z)$). A vector function $g : X \rightarrow Y^\bullet$ is said to be \textit{proper} whenever its \textit{domain} $\dom g = \{x \in X: g(x) \in Y\}$ is nonempty. For $\lambda \in K^*$ we denote by $(\lambda g) :X \rightarrow \overline{\R}$ the function defined by $(\lambda g)(x) =
\langle \lambda, g(x) \rangle$. The vector function $g$ is called \textit{$K$-convex} if for all $x,y \in X$ and all $t \in [0,1]$ one has $g(tx + (1-t)y) \leq_K t g(x) + (1-t) g(y)$. For $g : X \rightarrow Y^\bullet$ a $K$-convex function the set $Q=\{x \in X : g(x) \in -K\}$ turns out to be convex.

\section{Error bounds: the scalar case}

We say that the \textit{Slater qualification condition} holds for the inequality \eqref{ineq-sc} if
\begin{equation}\label{slater_sc}
\mbox{there exists} \ x_0 \in X \ \mbox{such that} \ f(x_0) < 0.
\end{equation}

In this section we show that the fulfillment of the Slater qualification condition combined with the boundedness of $S$ ensures the existence of global error bounds for \eqref{ineq-sc}. To this aim we give first an equivalent characterization of the existence of error bounds by means of conjugate functions.

\begin{lemma}\label{equiv-err-sc} Suppose that $f:X\rightarrow\B$ is a proper and convex function such that $S=\{x\in X:f(x)\leq 0\}$ is nonempty. Then the global error bound holds for the inequality \eqref{ineq-sc} with constant $\alpha > 0$ if and only if
\begin{equation}\label{ineq-equiv-err-sc}\min_{\lambda\in[0,1]}(\lambda f)^*(x^*)\leq\sigma_S(x^*) \ \forall x^*\in X^* \mbox{ with } \|x^*\|_{*}\leq 1/\alpha.\end{equation}\end{lemma}

\noindent {\bf Proof.} The key observation is that the inequality \eqref{err-b-sc} is fulfilled with constant $\alpha > 0$ if and only if \begin{equation}\label{ineq}(1/\alpha d(\cdot, S))^*\geq (f_+)^*.\end{equation} Indeed, the direct implication is trivial, while for the reverse one we use the fact that $f_+(x)\geq(f_+)^{**}(x)$ for all $x\in X$ and that $d(\cdot, S)$ is real-valued, convex and continuous, which allows to apply for it the Fenchel-Moreau Theorem. As $d(\cdot, S)=\|\cdot\|\Box\delta_S$, we have for all $x^* \in X^*$ that
\begin{equation}\label{conj1}(1/\alpha d(\cdot, S))^*(x^*)=1/\alpha\big(\|\cdot\|^*(\alpha x^*)+\delta_S^*(\alpha x^*)\big)=1/\alpha\delta_{B_*(0,1)}(\alpha x^*)+\sigma_S(x^*).\end{equation} Further, for all $x^* \in X^*$ we have (cf. \cite[Lemma 45.1]{simons}, see also \cite{bw}) \begin{equation}\label{conj2}(f_+)^*(x^*)=\min_{\lambda\in[0,1]}(\lambda f)^*(x^*).\end{equation} The result follows now from \eqref{ineq}, \eqref{conj1} and \eqref{conj2}.\hfill{$\Box$}

\begin{remark}\label{rem-equiv-sc}\rm (a) When $f$ is additionally lower semicontinuous one can alternatively use \cite[Theorem 2.1]{Zal-baikal} for proving the statement in Lemma \ref{equiv-err-sc}.

(b) For a fixed $x^*\in X^*$ consider the primal optimization problem
$$(P_{x^*}) \hspace{2cm} \inf_{x\in S}\langle -x^*,x\rangle$$
and its \textit{Lagrange dual problem}
$$(D_{x^*}) \hspace{2cm} \sup_{\lambda\geq 0}\inf_{x\in X}\{\langle -x^*,x\rangle+\lambda f(x)\}.$$ Since weak duality always holds, that is $v(P_{x^*})\geq v(D_{x^*})$, where $v(P_{x^*}), v(D_{x^*})$ are the optimal objective values of $(P_{x^*})$ and $(D_{x^*})$, respectively, one can easily derive the inequality
$$\sigma_S(x^*) = -v(P_{x^*}) \leq -v(D_{x^*}) = \inf_{\lambda\geq 0}(\lambda f)^*(x^*).$$
Hence \eqref{ineq-equiv-err-sc} in Lemma \ref{equiv-err-sc} can be equivalently written as
\begin{equation}\label{equivlem1}
\min_{\lambda\in[0,1]}(\lambda f)^*(x^*) = \sigma_S(x^*) \ \forall x^*\in X^* \mbox{ with } \|x^*\|_{*}\leq 1/\alpha.
\end{equation}
This means that the global error bound holds for \eqref{ineq-sc} with constant $\alpha > 0$ if and only if for all $x^* \in 1/\alpha B_*(0,1)$ one has $v(P_{x^*})=v(D_{x^*})$ and the dual $(D_{x^*})$ has an optimal solution $\bar \lambda$ in the interval $[0,1]$ (this can be seen as a \textit{sharp strong duality} statement for the primal-dual pair $(P_{x^*})-(D_{x^*})$).

(c) One can easily notice that when $f$ is proper and convex and the Slater qualification condition for \eqref{ineq-sc} is fulfilled, then for $(P_{x^*})-(D_{x^*})$ strong duality holds for all $x^*\in X^*$ (see, for instance, \cite{b-hab, bgw-book, Zal-carte}), which is nothing else than
\begin{equation}\label{eq-sigma-slater}\sigma_S(x^*) = -v(P_{x^*}) = -v(D_{x^*}) = \min_{\lambda\geq 0}(\lambda f)^*(x^*).
\end{equation}
\end{remark}

We come now to the first global error bound result for \eqref{ineq-sc}, for the proof of which we use conjugate duality techniques, but also an useful characterization of the continuity of the conjugate of a function given by Rockafellar in \cite{rock-level-cont}.

\begin{theorem}\label{err-sc-lsc} Suppose that $f:X\rightarrow\B$ is a proper, convex and lower semicontinuous function such that the Slater qualification condition
for \eqref{ineq-sc} is fulfilled and $S=\{x\in X:f(x)\leq 0\}$ is a bounded set. Then for the inequality \eqref{ineq-sc} the global error bound holds.\end{theorem}

\noindent {\bf Proof.} We show first that for all $r \in \R$ the lower level set $L(f,r)$ is bounded. To this aim we fix $r \in \R$. Let $x_0\in X$ be such that $f(x_0)<0$ and $M\geq 0$ fulfilling $\|x\|\leq M$ for all $x\in S$. There exists a sufficiently small $\lambda\in (0,1)$ such that the inequality $f(x_0)+\lambda(r-f(x_0))<0$ is fulfilled. Take now an arbitrary element $x\in L(f,r)$. The function $f$ being convex we get $$f(x_0+\lambda(x-x_0))\leq(1-\lambda)f(x_0)+\lambda f(x)\leq f(x_0)+\lambda(r-f(x_0))<0,$$ hence $\|x_0+\lambda(x-x_0)\|\leq M$, which ensures that $\|x\|\leq 1/\lambda(\|x_0\|+M)+\|x_0\|$. Therefore $L(f,r)$ is bounded. As $r \in \R$ was arbitrarily chosen, in virtue of \cite[Theorem 7A(a)]{rock-level-cont}, the conjugate $f^*$ is finite and strongly continuous at $0$. Moreover, the Slater qualification condition ensures that $f^*(0)\geq -f(x_0) > 0$.

Suppose in the following that the global error bound does not hold for the inequality \eqref{ineq-sc}. Applying Lemma \ref{equiv-err-sc} and taking into consideration Remark \ref{rem-equiv-sc}(b) and (c), it follows that for all $\alpha>0$ there exist $x_{\alpha}^*\in X^*$, $\|x_{\alpha}^*\|_{*}\leq 1/\alpha$ and $r_{\alpha}\in\R$ such that $$\min_{\lambda\in[0,1]}(\lambda f)^*(x_{\alpha}^*)>r_{\alpha}>\min_{\lambda\geq 0}(\lambda f)^*(x_{\alpha}^*).$$ By taking $\alpha:=n$ ($n\in\N$), we obtain the existence of sequences $x_n^*\in X^*$, $\|x_n^*\|_{*}\leq 1/n$, $r_n\in\R$ and $\lambda_n\in\R$, $\lambda_n>1$ ($n \in \N$) such that \begin{equation}\label{rn1}r_n>(\lambda_n f)^*(x_n^*)=\lambda_n f^*(1/{\lambda_n}x_n^*) \ \forall n\in\N\end{equation} and \begin{equation}\label{rn2}r_n<(\lambda f)^*(x_n^*)=\lambda f^*(1/{\lambda}x_n^*)  \ \forall n\in\N \ \forall \lambda\in(0,1].\end{equation} In the argumentation below we use that $f^*$ is finite and continuous at $0$, $f^*(0)>0$ and $x_n^*\rightarrow 0$.

We consider two cases: the first one when the sequence $\lambda_n$ is unbounded. Then there exists a subsequence $\lambda_{n_k}$ ($k\in\N$) such that $\lambda_{n_k}\rightarrow\infty$ ($k\rightarrow\infty$). From \eqref{rn1} we get $r_{n_k}\rightarrow\infty$ ($k\rightarrow\infty$), which contradicts \eqref{rn2}.

Suppose now that the sequence $\lambda_n$ is bounded. There exists a convergent subsequence $\lambda_{n_i}$ ($i\in\N$) such that $\lambda_{n_i}\rightarrow \overline{\lambda}\in[1,\infty)$ ($i\rightarrow\infty$). From \eqref{rn1} and \eqref{rn2} we obtain $$\overline{\lambda} f^*(0)\leq\liminf_{i\rightarrow\infty}r_{n_i}\leq \limsup_{i\rightarrow\infty}r_{n_i}\leq \lambda f^*(0) \ \forall\lambda\in(0,1],$$ which is, of course, impossible.

Thus our statement, that the global error bound does not hold for the inequality \eqref{ineq-sc}, is false and the proof is complete.\hfill{$\Box$}

We show in the following that in the above theorem the lower semicontinuity of the function $f$ can be dropped. To this aim we work with the lower semicontinuous hull of $f$. For other considerations concerning the relation between the existence of global error bounds for \eqref{ineq-sc} and the existence of global error bounds for a similar inequality, where $f$ is replaced by $\cl f$, we refer to \cite{hu}.

\begin{theorem}\label{err-sc} Suppose that $f:X\rightarrow\B$ is a proper and convex function such that the Slater qualification condition
for \eqref{ineq-sc} is fulfilled and $S=\{x\in X:f(x)\leq 0\}$ is a bounded set. Then for the inequality \eqref{ineq-sc} the global error bound holds.\end{theorem}

\noindent {\bf Proof.} Let $x_0\in X$ be such that $f(x_0)<0$. We show first that the Slater qualification condition guarantees the following equality  \begin{equation}\label{eq-cl-f}\cl S = \{x \in X : \cl f(x) \leq 0\}.
\end{equation}
Since the inclusion ``$\subseteq$'' is obvious, we prove only the reverse one. Let $x\in X$ be such that $\cl f(x)\leq 0$. This means that $(x,0)\in\epi(\cl f)=\cl\big(\epi f\big)$. Hence there exist sequences $x_n\in X$, $r_n\in\R$ ($n\in\N$) such that $f(x_n)\leq r_n$ for all $n\in\N$ and $(x_n,r_n)\rightarrow(x,0)$ ($n\rightarrow\infty$). We can suppose without losing the generality that $r_n\leq 1/n^2$ for all $n \in \N$. Since $f(x_0)<0$, there exists $n_0\in \N$ such that $f(x_0)+(1-1/n)1/n<0$ for all $n\geq n_0$. Define the sequence $y_n:=(1/n)x_0+(1-1/n)x_n$ ($n\in\N$). The convexity of the function $f$ ensures $$f(y_n)\leq (1/n)f(x_0)+(1-1/n)f(x_n)\leq (1/n)f(x_0)+(1-1/n)1/n^2<0 \ \forall n\geq n_0,$$ hence $y_n \in S$ for all $n\geq n_0$. Since $y_n\rightarrow x$ ($n\rightarrow\infty$), we conclude that \eqref{eq-cl-f} holds.

We prove that the global error bound holds for the inequality \begin{equation}\label{ineq-sc-cl}\cl f(x)\leq 0, \ x\in X\end{equation} that is there exists a constant $\alpha>0$ such that \begin{equation}\label{err-b-sc-cl}d(x,\{y \in X: \cl f(y) \leq 0\})\leq \alpha [\cl f(x)]_+ \ \forall x\in X.\end{equation} We consider to this aim two cases: the first one when $\cl f$ is not proper. Due to \cite[Proposition 2.4]{EkTem} we get $\cl f(x)=-\infty$ for $x\in\dom(\cl f)$ and $\cl f(x)=+\infty$ for $x\not\in\dom(\cl f)$. Then $\{y \in X: \cl f(y) \leq 0\} = \dom(\cl f)$ and thus \eqref{err-b-sc-cl} holds for arbitrary $\alpha>0$.

In case the function $\cl f$ is proper, relation \eqref{eq-cl-f} and the fact that $\cl f(x_0) \leq f(x_0) < 0$ guarantee that Theorem \ref{err-sc-lsc} can be applied for the function $\cl f$, thus the global error bounds holds for the inequality \eqref{ineq-sc-cl}. This means that there exists $\alpha>0$ such that \eqref{err-b-sc-cl} holds. Since for all $x \in X$ one has  $f(x)\geq\cl f(x)$ and (cf. \eqref{eq-cl-f}) $d(x,\{y \in X: \cl f(y) \leq 0\})=d(x,\cl S)=d(x,S)$, the global error bound holds for the inequality \eqref{ineq-sc} with the same constant $\alpha > 0$. \hfill{$\Box$}

In general the Slater qualification condition is not enough in order to guarantee the existence of global error bounds (see \cite[Example 2]{lew-pang}). However, for particular convex inequality systems one can renounce to the boundedness condition. We prove in the following that in the above theorem the assumption, that the lower level set of $f$ at level $0$ is bounded, can be removed in case $f:\R^m\rightarrow\R$ ($m\in\N$) is a convex quadratic function. Here we consider $\R^m$ to be endowed with an arbitrary norm.

\begin{theorem}\label{err-sc-quadr} Let $A$ be a $m\times m$ symmetric positive semidefinite matrix ($m\in \N$), $b \in \R^m$, $c \in \R$ and for the function $f:\R^m\rightarrow\R$, defined by $f(x)=1/2\langle x,Ax\rangle+\langle b,x\rangle-c$, let us assume that the Slater qualification condition
for \eqref{ineq-sc} is fulfilled. Then for the inequality \eqref{ineq-sc} the global error bound holds.\end{theorem}

\noindent {\bf Proof.} As in the proof of Theorem \ref{err-sc-lsc}, we suppose that the global error bound does not hold for the inequality \eqref{ineq-sc}, hence there exist sequences $x_n^*\in \R^m$, $\|x_n^*\|\leq 1/n$, $r_n\in\R$, $\lambda_n\in\R$, $\lambda_n>1$ ($n \in \N$) such that \begin{equation}\label{rn1'}r_n>(\lambda_n f)^*(x_n^*)=\lambda_n f^*(1/{\lambda_n}x_n^*) \ \forall n\in\N\end{equation} and \begin{equation}\label{rn2'}r_n<(\lambda f)^*(x_n^*)=\lambda f^*(1/{\lambda}x_n^*) \ \forall n\in\N \ \forall \lambda\in(0,1].\end{equation} The Slater condition ensures that $f^*(0)>0$. Moreover, as $\dom f^*=b+A(\R^m)$ (cf. \cite[Chapter X, Example 1.1.4]{hiriart-lem1}), from \eqref{rn1'} we get $(1/{\lambda_n})x_n^*\in b+A(\R^m)$ for all $n\in\N$. Since $(1/{\lambda_n})x_n^*\rightarrow 0$ ($n\rightarrow+\infty$), we obtain
$$0\in\cl(b+A(\R^m))=b+A(\R^m)=\dom f^*.$$
This implies that $\dom f^*=b+A(\R^m)=A(\R^m)$ and $f^*(0)\in\R$. Thus $\ri \dom f^* = \ri A(\R^m)=A(\R^m)=\dom f^*$ (cf. \cite[Theorem 6.6]{Rock-carte}). From $(1/{\lambda_n})x_n^*\in b+A(\R^m)=A(\R^m)$ we derive $x_n^*\in A(\R^m)$ and  so $(1/\lambda)x_n^*\in A(\R^m)$ for all $n\in\N$ and all $\lambda\in(0,1]$. Using the fact that $f^*$ is continuous relative to $\ri \dom f^*= A(\R^m)$ (cf. \cite[Theorem 10.1]{Rock-carte}), the proof can be continued in the lines of the second part of the proof of Theorem \ref{err-sc-lsc}.\hfill{$\Box$}

\begin{remark}\rm Luo and Luo proved in \cite{luo}, by using some results from the linear algebra, the existence of global error bound results also for inequality systems of $k$ ($k \in \N$) convex quadratic functions. At this moment we are not aware of how the techniques used in the proof of Theorem \ref{err-sc-quadr} can be extended to this more general situation.
\end{remark}

\section{Error bounds: from the scalar to the vector case}

In this section we consider a further real normed space $(Y,\|\cdot\|)$, partially ordered by a \textit{convex closed cone} $K \subseteq Y$ having a \textit{nonempty interior}, and $g : X \rightarrow Y^\bullet$ a proper vector function such that $Q=\{x\in X:g(x)\in -K\}$ is nonempty. We will provide some existence results for the global error bound of the inequality system \eqref{ineq-vect}, which we deduce from the scalar case investigated above.

To this aim we make use of the \textit{oriented distance function}, which is a special scalarization function introduced by Hiriart-Urruty in \cite{hir-bull, hir-mor}. For $A\subseteq Y$ this function is defined by
$$\Delta_A:Y\rightarrow\B, \Delta_A(y)=d(y,A)-d(y,Y\setminus A) \ \mbox{for all} \ y\in Y.$$
Let us recall in the following the properties of the oriented distance function which we will be used throughout this section (see \cite[Proposition 3.2]{zaf}). Suppose that $A$ is a nonempty, convex  and closed set such that $A \neq Y$. Then $\Delta_A$ is real-valued, convex, while
$$\{y\in Y:\Delta_A(y)\leq 0\}=A \ \mbox {and} \ \inte A\subseteq\{y\in Y:\Delta_A(y)<0\}.$$
If, additionally, $A$ is a cone, then $\Delta_A$ is $K$-decreasing.

We say that the \textit{Slater qualification condition} holds for the inequality \eqref{ineq-vect} if
\begin{equation}\label{slater_vect}
\mbox{there exists} \ x_0 \in X \ \mbox{such that} \ g(x_0) \in -\inte K.
\end{equation}

The following result was first proved by Robinson in \cite{robin} (in case the function $g$ is defined on a nonempty convex subset of $X$). We give here an alternative proof for it, which relies on Theorem \ref{err-sc} and makes use of the oriented distance function.

\begin{theorem}\label{err-vect} Suppose that $g:X\rightarrow Y^\bullet$ is a proper and $K$-convex function such that the Slater qualification condition
for \eqref{ineq-vect} is fulfilled and $Q=\{x\in X:g(x)\in -K\}$ is a bounded set. Then for the inequality system \eqref{ineq-vect} the global error bound holds.\end{theorem}

\noindent {\bf Proof.} When $K=Y$, then \eqref{err-b-vect} holds for an arbitrary $\alpha > 0$. Assume in the following that $K \neq Y$.

Consider the function $f:X\rightarrow \overline{\R}$ defined by $f(x)=\Delta_{-K}(g(x))$ for all $x\in X$. The properties of the oriented distance function guarantee that $f$ is a proper and convex function. Moreover, \begin{equation}\label{f-g-vect1}\{x\in X:f(x)\leq 0\}=\{x\in X:\Delta_{-K}(g(x))\leq 0\}=\{x\in X:g(x)\in -K\}=Q,\end{equation}
while the Slater qualification condition guarantees that $f(x_0) = \Delta_{-K}(g(x_0)) < 0$. This means that all the hypotheses of Theorem \ref{err-sc} are verified, hence there exists $\alpha>0$ such that
\begin{equation}\label{f-g-vect2}
d(x,Q) = d(x,\{y \in X : f(y) \leq 0\}) \leq \alpha f(x)_+ \ \forall x \in X.
\end{equation}
We close the proof by showing that the global error bound holds for the inequality system \eqref{ineq-vect} with the same constant $\alpha$. Take an arbitrary $x\in X$. If $g(x) = \infty_K$, then \eqref{err-b-vect} is obviously fulfilled. Further assume that $g(x) \in Y$. If $x\in Q$, that is, $g(x)\in-K$, then obviously $d(x,Q) = 0 = \alpha d(g(x),-K)$. If $x\not\in Q$, that is, $g(x)\in Y \setminus (-K)$, then $f(x)>0$ (cf. \eqref{f-g-vect1}). From \eqref{f-g-vect2} we get $d(x,Q) \leq \alpha f(x)$ and the conclusion follows, since in this case $f(x)=d(g(x),-K)$.\hfill{$\Box$}

\section{Sharpening the error bound result of Robinson}

In this section we work in the setting of the previous section and give an alternative proof for the existence of global error bounds for \eqref{ineq-vect}, succeeding meanwhile to sharpen the statement of Robinson in \cite{robin} concerning the bound $\alpha > 0$. Recall that, under the assumption that the Slater qualification condition for \eqref{ineq-vect} is fulfilled at $x_0 \in X$ and that $Q$ is bounded, Robinson proved that \eqref{err-b-vect} is fulfilled for $\alpha =\diam Q/\delta $, where $\delta > 0$ is such that
$\delta B(0,1) \subseteq g(x_0)+K$ and $\diam Q := \sup\{\|y-z\| : y,z \in Q\}$ is the \textit{diameter} of the set $Q$.

According to \cite{robin}, when $x \in X$ is such that $g(x) \in Y \setminus (-K)$, then for  $\rho:=d(g(x),-K)$ $>0$ and $\lambda:=\rho/(\rho+\delta)\in (0,1)$, one has $(1-\lambda)x+\lambda x_0\in Q$, which means that the set $Q$ does not reduce to a singleton, that is $\diam Q>0$.

We start as in the scalar case with an equivalent characterization of the existence of error bounds by means of conjugate functions.

\begin{lemma}\label{equiv-err-vect} Suppose that $g:X\rightarrow Y^\bullet$ is a proper and $K$-convex function such that $Q=\{x\in X:g(x)\in -K\}$ is nonempty. Then the global error bound holds for the inequality \eqref{ineq-vect} with constant $\alpha > 0$ if and only if
\begin{equation}\label{ineq-equiv-err-vect}\min_{\substack{\lambda\in K^*\\ \|\lambda\|_{*}\leq 1}}(\lambda g)^*(x^*)\leq\sigma_Q(x^*) \ \forall x^*\in X^* \mbox{ with } \|x^*\|_{*}\leq 1/\alpha.\end{equation}\end{lemma}

\noindent {\bf Proof.} Relation \eqref{err-b-vect} is equivalent to \begin{equation}\label{ineq-vect-dem}(1/\alpha d(\cdot, Q))^*\geq f^*,\end{equation} where $f:X\rightarrow \overline{\R}$, $f=d(\cdot,-K)\circ g$. One can easily show that the function $d(\cdot, -K)$ is $K$-increasing, hence $f$ is proper and convex. Moreover, since $d(\cdot,-K)$ is continuous, we can apply \cite[Theorem 3.5.2(a)]{bgw-book} in order to compute the conjugate of $f$. For all $x^*\in X^*$ we get \begin{equation}\label{fstar}f^*(x^*)=\min_{\lambda\in K^*}\Big[\big(d(\cdot,-K)\big)^*(\lambda)+(\lambda g)^*(x^*)\Big].\end{equation} Since $d(\cdot,-K)=\|\cdot\|\Box\delta_{-K}$, we get for all $\lambda\in K^*$
$$\big(d(\cdot,-K)\big)^*(\lambda)=(\|\cdot\|)^*(\lambda)+\sigma_{-K}(\lambda),$$
which is equal to $0$, for $\|\lambda\|_{*}\leq 1$, being $+\infty$, otherwise. Hence, \begin{equation}\label{fstar2}f^*(x^*)=\min_{\substack{\lambda\in K^*\\\|\lambda\|_{*}\leq 1}}(\lambda g)^*(x^*) \ \forall x^*\in X^*.\end{equation} As $$(1/\alpha d(\cdot, Q))^*=1/\alpha\big(\|\cdot\|^*(\alpha x^*)+\delta_Q^*(\alpha x^*)\big)=1/\alpha\delta_{B_*(0,1)}(\alpha x^*)+\sigma_Q(x^*),$$ the result follows from \eqref{ineq-vect-dem} and \eqref{fstar2}.\hfill{$\Box$}

\begin{remark}\label{rem-equiv-vect}\rm (a) One can notice that the equivalence in the above lemma remains true even if $K$ fails to be closed or to have a nonempty interior.

(b) For a fixed $x^*\in X^*$ consider the primal optimization problem

$$(P^v_{x^*}) \ \ \ \inf_{x\in Q}\langle -x^*,x\rangle$$ and its \textit{Lagrange dual problem} $$(D^v_{x^*}) \ \ \ \sup_{\lambda\in K^*}\inf_{x\in X}\{\langle -x^*,x\rangle+(\lambda g)(x)\}.$$ Since weak duality always holds, that is $v(P^v_{x^*})\geq v(D^v_{x^*})$, where $v(P^v_{x^*}), v(D^v_{x^*})$ are the optimal objective values of $(P^v_{x^*})$ respectively $(D^v_{x^*})$, one can easily derive the inequality
$$\sigma_Q(x^*) = -v(P^v_{x^*}) \leq -v(D^v_{x^*}) = \inf_{\lambda\in K^*}(\lambda g)^*(x^*).$$
Hence \eqref{ineq-equiv-err-vect} in Lemma \ref{equiv-err-vect} can be equivalently written as
\begin{equation}\label{equivlem6}
\min_{\substack{\lambda\in K^*\\ \|\lambda\|_{*}\leq 1}}(\lambda g)^*(x^*) = \sigma_Q(x^*) \ \forall x^*\in X^* \mbox{ with } \|x^*\|_{*}\leq 1/\alpha.
\end{equation}
This means that the global error bound holds for \eqref{ineq-vect} with constant $\alpha > 0$ if and only if for all $x^* \in 1/\alpha B_*(0,1)$ one has $v(P^v_{x^*})=v(D^v_{x^*})$ and the dual $(D^v_{x^*})$ has an optimal solution $\bar \lambda$ in the set $K^* \cap B_*(0,1)$ (this can be seen as a \textit{sharp strong duality} statement for the primal-dual pair $(P^v_{x^*})-(D^v_{x^*})$).

(c) One can easily notice that when $g$ is proper and $K$-convex and the Slater qualification condition for \eqref{ineq-vect} is fulfilled, then for $(P^v_{x^*})-(D^v_{x^*})$ strong duality holds for all $x^*\in X^*$ (see, for instance, \cite{b-hab, bgw-book, Zal-carte}), which is nothing else than
\begin{equation}\label{eq-sigma-slater-vect}\sigma_Q(x^*) = -v(P^v_{x^*}) = -v(D^v_{x^*}) = \min_{\lambda\in K^*}(\lambda g)^*(x^*).
\end{equation}
\end{remark}

In the proof of the following statement we use a \textit{sharp Lagrange multiplier} result due to Simons.

\begin{theorem}\label{err-vect_1}Suppose that $g:X\rightarrow Y^\bullet$ is a proper and $K$-convex function such that the Slater qualification condition
for \eqref{ineq-vect} is fulfilled at $x_0 \in X$, i.e. $g(x_0)\in-\inte K$, and $Q=\{x\in X:g(x)\in -K\}$ is a bounded set. Then for the inequality system \eqref{ineq-vect} the global error bound holds with $$\alpha=\frac{\diam Q}{d(g(x_0),Y\setminus ({-K}))}.$$\end{theorem}

\noindent {\bf Proof.} Taking into account Remark \ref{rem-equiv-vect}(b), it is enough to show that for all $x^*\in X^*$ with $\|x^*\|_{*}\leq 1/\alpha$, strong duality holds for the primal-dual pair $(P^v_{x^*})-(D^v_{x^*})$ and that $(D^v_{x^*})$ has an optimal solution $\lambda\in K^*$ with $\|\lambda\|_{*}\leq 1$.

Take an arbitrary $x^*\in X^*$ with $\|x^*\|_{*}\leq 1/\alpha$. Since the Slater qualification condition is fulfilled, we can apply \cite[Theorem 6.6]{simons}. It follows that strong duality holds for the pair $(P^v_{x^*})-(D^v_{x^*})$ and $(D^v_{x^*})$ has an optimal solution $\lambda\in K^*$ with $$\|\lambda\|_{*}\leq \inf_{\substack {x\in X\\g(x)\in-\inte K}}\frac{\langle-x^*,x\rangle-\inf_{u\in Q}\langle-x^*,u\rangle}{d(g(x),Y\setminus {-K})}$$$$\leq\frac{\sup_{u\in Q}\langle x^*,u\rangle-\langle x^*,x_0\rangle}{d(g(x_0),Y\setminus {-K})}=\frac{\sup_{u\in Q}\langle x^*,u-x_0\rangle}{d(g(x_0),Y\setminus {-K})}\leq \frac{(1/\alpha)\diam Q}{d(g(x_0),Y\setminus {-K})}=1,$$ and the proof is complete.\hfill{$\Box$}

\begin{remark}\rm For $x_0 \in X$ with $g(x_0) \in -\inte K$ we proved that for \eqref{ineq-vect} the global error bound holds with $\alpha_{BC}:=\diam Q/ d(g(x_0),Y\setminus ({-K}))$, while Robinson in \cite{robin} got as a bound for the same inequality system $\alpha_R(\delta) :=\diam Q/\delta $, where $\delta>0$ is such that $\delta B(0,1) \subseteq g(x_0)+K$. In the following we prove that
$$\alpha_{BC} = \inf\{\alpha_R(\delta) : \delta > 0, \delta B(0,1) \subseteq g(x_0)+K\},$$
which actually means proving that
\begin{equation}\label{eqfinal} 
d(g(x_0),Y\setminus ({-K})) = \sup\{\delta > 0 : \delta B(0,1) \subseteq g(x_0)+K\}.
\end{equation}
Take first an arbitrary $\delta>0$ such that $\delta B(0,1) \subseteq g(x_0)+K$. Then $d(g(x_0),Y\setminus (-K))\geq \delta$. Indeed, if there exists $y_0\in Y\setminus (-K)$ such that $\|g(x_0)-y_0\|<\delta$, then $g(x_0)-y_0\in\delta B(0,1) \subseteq g(x_0)+K$, hence $y_0\in -K$, which is a contradiction. Thus
$$d(g(x_0),Y\setminus ({-K})) \geq \delta \ \mbox{for all} \ \delta > 0 \ \mbox{with} \ \delta B(0,1) \subseteq g(x_0)+K.$$
Take now an arbitrary $\delta > 0$ such that $\delta < d(g(x_0),Y\setminus (-K))$. Then one has that $\delta B(0,1) \subseteq g(x_0)+K$. Indeed, if there exists
$y_0\in Y\setminus (-K)$ such that $\|y_0-g(x_0)\|\leq\delta$, then $\|y_0-g(x_0)\| < d(g(x_0),Y\setminus {-K})\leq \|y_0-g(x_0)\|$, which is a contradiction.
Thus \eqref{eqfinal} holds and $\alpha_{BC}$ proves to be the infimum over the family of bounds proposed by Robinson in \cite{robin}.
\end{remark}

\section{Conclusion and further research}

We have shown that the theory of conjugate duality can be successfully applied in order to get existence results concerning global error bounds for convex inequality systems. We investigated in the first part the scalar case and then we have proposed a bridge between the scalar and the vector case via the \textit{oriented distance function} introduced by Hiriart-Urruty. In the last section we computed by means of a Lagrange multiplier result due to Simons a bound which sharpens the ones given by Robinson in the context of error bounds for convex inequality systems defined by vector functions.

An interesting future research topic in this area could be to find out if the conjugate duality techniques used in this paper can be implemented in case of error bounds defined by multifunctions. More precisely, having $\Gamma :X\rightrightarrows Y$ a multifunction, where $X,Y$ are real normed spaces, we say that $\Gamma$ has a \textit{global error bound at $x_0 \in \dom\Gamma$} provided that there exists $\alpha>0$ such that $$d(y,\Gamma(x_0))\leq\alpha d(x_0,\Gamma^{-1}(y)) \ \forall y\in Y.$$ This is a generalization of the notions considered in this paper. We refer to \cite{li-sing, zal-mor, klatte} for conditions which guarantee the existence of error bounds in this context.

Another direction, which could be of interest, is to analyze if instead of the Slater qualification condition, which requires the nonemptiness of the interior of the cone $K$, some weaker conditions could be considered, in order to guarantee the existence of error bounds. Recall that there exist generalizations of the classical interior, like the \emph{algebraic interior}, the \emph{strong quasi-relative interior} and the \emph{quasi-relative interior}, which play an important role in the formulation of regularity conditions ensuring strong duality in convex optimization; see \cite{b-hab, bgw-book, cs, Zal-carte}.

Finally, it could be challenging to see if the technique used in the proof of Theorem \ref{err-sc-quadr} can be generalized to $k$ convex quadratic functions (with the corresponding global error bound notion in the vector case). We know that a similar result remains valid in this case, too (cf. \cite[Theorem 3.1]{luo}).

\end{document}